\documentclass[aos,MSNbibl,seceqn,dvips]{arximspdf}

\doi{10.1214/13-AOS1179} 
\volume{42}
\issue{3}
\pubyear{2014}
\firstpage{1131}
\lastpage{1144}

\makeatletter
\renewcommand{\mid}{|}
\newcommand{\wC}{\widehat{C}}
\newcommand{\Weta}{\widetilde{\eta}}
\newcommand{\Wphi}{\widetilde{\phi}}
\newcommand{\Beta}{\overline{\eta}}
\newcommand{\tols}{\stackrel{\mathcal{L}-\mathrm{s}}{\longrightarrow}}
\newtheorem{prop}{Proposition}[section]
\newtheorem{theo}[prop]{Theorem}
\newproclaim{rem}[prop]{Remark}
\newproclaim{ass}{Assumption}
\makeatother

\begin{document}
\begin{frontmatter}

\title{A remark on the rates of convergence for integrated volatility
estimation in\\ the presence of jumps}
\runtitle{Rates for integrated volatility}

\begin{aug}
\author[a]{\fnms{Jean}~\snm{Jacod}\corref{}\thanksref{t1}\ead[label=e1]{jean.jacod@upmc.fr}}
\and
\author[b]{\fnms{Markus}~\snm{Reiss}\ead[label=e2]{mreiss@math.hu-berlin.de}}
\runauthor{J. Jacod and M. Reiss}
\affiliation{UPMC (Universit\'e Paris-6) and Humboldt-Universit\"at zu Berlin}
\address[a]{Institut de Math\'ematiques de Jussieu\\
UPMC (Universit\'e Paris-6)\\
4 Place Jussieu\\
75005-Paris\\
France\\
\printead{e1}} 
\address[b]{Institut f\"ur Mathematik\\
Humboldt-Universit\"at zu Berlin\\
Unterden Linden, 6\\
10099-Berlin\\
Germany\\
\printead{e2}}
\end{aug}
\thankstext{t1}{Supported in part by a Humboldt Research Award.}

\received{\smonth{9} \syear{2012}}
\revised{\smonth{9} \syear{2013}}

%
\begin{abstract}
The optimal rate of convergence of estimators of the integrated
volatility, for a discontinuous It\^o semimartingale sampled at regularly
spaced times and over a fixed time interval, has been a long-standing
problem, at least when the jumps are not summable. In this paper, we
study this optimal rate, in the minimax sense and for appropriate
``bounded'' nonparametric classes of semimartingales. We show that,
if the $r$th powers of the jumps are summable for some $r\in[0,2)$,
the minimax rate is equal to $\min(\sqrt{n},(n\log n)^{(2-r)/2})$, where
$n$ is the number of observations.
\end{abstract}

%
\end{frontmatter}

\section{Introduction}\label{sec-Intro}

Let $X$ be a one-dimensional It\^o semimartingale, which in particular
means that its ``continuous martingale part'' has the form
\[
X^c_t=\int_0^t
\sigma_s \,dW_s,
\]
where $W$ is a standard Brownian motion, and the process $\sigma_t$ is
optional and (locally) squared integrable.

One of the long-standing problems is the estimation of the so-called
integrated volatility, say at time $1$, that is of the variable
$C_1=\int_0^1c_s \,ds$, where $c_t=\sigma_t^2$ is the (squared) volatility,
on the basis of discrete observations of $X$. A huge number of papers have
been devoted to this question already, in various situations: when the
process is continuous (so $X$ is the sum of $X^c$ above, plus possibly a
drift term), or when it has jumps; when
the process $X$ is ``perfectly'' observed, or contaminated by noise;
when the sampling times are equi-spaced, or when they are
irregularly spaced.

Below, we focus on the basic case, where the sampling is at
regularly spaced times $i/n$ for $i=0,\ldots,n$, and when $X_{i/n}$
is observed without noise. Even in this simple situation, the question
of the
``optimal'' rate of convergence of estimators toward $C_1$, as $n\to
\infty$,
is unanswered so far, when there are jumps which are ``too active.''

More precisely, estimators are known, which converge to $C_1$ with the rate
$\sqrt{n}$, in the continuous case (the realized volatility, or
``approximate quadratic variation'' at time $1$, achieves this rate), and
also when $X$ has jumps with a degree of activity, or
Blumenthal--Getoor index,
less than $1$. This rate is optimal (in a minimax sense), for the following
reason: if $X=\sigma W$ where $c=\sigma^2$ is a constant, so $C_1=c$,
we have a
purely parametric setting for which the local asymptotic normality (LAN)
holds with rate $\sqrt{n}$, and the realized volatility is indeed the
MLE in
this case.

However,\vspace*{1pt} when the degree $r$ of jump activity is larger than $1$, the best
rates found in the literature are of the form $n^{((2-r)/2)-\varepsilon}$ for
$\varepsilon>0$ arbitrarily small (see below for more details). The difficulty
comes of course from the essentially nonparametric feature of the
problem, since we do not want to specify the law of the process $X$, apart
from the fact that it is an It\^o semimartingale, plus possibly some
boundedness assumptions on its characteristics. In a purely
parametric problem, for example, when $X$ is a L\'evy process with a
\textit{known} L\'evy measure and the only unknown parameters are the variance
$c$ of the Gaussian part, and possibly the drift, then again the rate
$\sqrt{n}$ is available for estimating $c$ (this rate is achieved by
the MLE,
under very general circumstances). There has been a considerable
interest in
providing also nonparametric estimators that converge at rate $\sqrt{n}$,
but as we
show here, this is in general impossible.

In this paper, a bound for the minimax rate is determined,
when the degree of activity is $r$ or smaller [the precise
definition of $r$ is given in Assumption~\ref{assLr} below, and is slightly
different from the usual Blumenthal--Getoor index]. We will see that
the best
possible rate is $(n\log n)^{(2-r)/2}$ when $r>1$ (and of course
$\sqrt{n}$ when $r\leq1$). It is interesting to notice that the truncated
realized volatility, which achieves the rate $n^{((2-r)/2)-\varepsilon}$ for
any prespecified $\varepsilon>0$ is indeed ``almost'' rate-optimal.

The paper is organized as follows: in Section~\ref{sec-K}, we state the
assumptions and review some known results. The results of this paper are
presented in Section~\ref{sec-R}, and the proofs are given in the last section.

\section{Some known results}\label{sec-K}

We consider a one-dimensional It\^o semimartingale $X$ on a filtered space
$(\Omega,\mathcal{F},(\mathcal{F}_t)_{t\geq0},\mathbb{P})$, which
is observed at regularly spaced times $\frac{i}n$ for
$i=0,1,\ldots,n$, over the (fixed) finite interval $[0,1]$. The characteristics
$(B,C,\nu)$ where $B$ is the drift, $C$ the integrated volatility and
$\nu$ the L\'evy system of~$X$ (see, e.g., Chapter~1 of \cite{JP}),
thus have the form
%
\begin{equation}
\label{1} B_t=\int_0^tb_s
\,ds,\qquad C_t=\int_0^tc_s
\,ds,\qquad \nu(dt,dx)=dt F_t(dx).
\end{equation}
Here, $b_t$ and $c_t$ are optional (or, predictable) processes, with
$c_t\geq0$,
and $F_t=F_{\omega,t}(dx)$ is an optional random measure, also called the
L\'evy measure, which accounts for the jumps of the process.

When $X$ is continuous, the canonical way for estimating $C_1$ is to
use the realized volatility, or approximate quadratic variation at time $1$:
%
\begin{equation}
\label{2} [X,X]^n_1=\sum
_{i=1}^n\bigl(\Delta^n_iX
\bigr)^2\qquad\mbox{where } \Delta^n_iX=X_{i/n}-X_{(i-1)/n},
\end{equation}
which converges in probability to $C_1$.
When further $\int_0^1b^2_s \,ds$ and $\int_0^1c_s^2 \,ds$ are a.s. finite,
we have the stable convergence in law at rate $\sqrt{n}$
%
\begin{equation}
\label{3} \sqrt{n} \bigl([X,X]^n_1-C_1
\bigr) \tols \mathcal{U}\qquad\mbox{where } \mathcal{U}=\sqrt{2}\int
_0^1c_s \,dW'_s,
\end{equation}
and where $W'$ is a standard Brownian motion, defined on an extension of
$(\Omega,\mathcal{F},(\mathcal{F}_t)_{t\geq0},\mathbb{P})$, and
which is independent of the $\sigma$-field $\mathcal{F}$: see, for example,
Theorem 5.4.2 in \cite{JP}.

When $X$ has jumps, the variables $[X,X]^n_1$ no longer converge to
$C_1$, but
to the ``full'' quadratic variation $[X,X]_1=
C_1+\sum_{s\leq1}(\Delta X_s)^2$, where $\Delta X_s=X_s-X_{s-}$
denotes the
jump size at time $s$. However, there are two known methods
to consistently estimate $C_1$:
\begin{longlist}[(2)]
\item[(1)] \textit{Truncated realized volatility}.  One chooses a sequence $v_n$
of positive truncation levels, typically of the form $v_n\asymp
1/n^\varpi$
for some $\varpi\in(0,1/2)$, and considers
%
\begin{equation}
\label{4} \wC(v_n)_1=\sum
_{i=1}^n\bigl(\Delta^n_iX
\bigr)^2 1_{\{|\Delta^n_iX|\leq v_n\}}.
\end{equation}

\item[(2)] \textit{Multipower variations}. One chooses an integer $k\geq2$, and
considers
%
\begin{equation}
\label{5} \wC(k,n)_1=\frac{1}{m_{2/k}^k} \sum
_{i=1}^{n-k+1} \prod_{j=0}^{k-1}\bigl|
\Delta^n_{i+j}X\bigr|^{2/k},
\end{equation}
where $m_p=\mathbb{E}(|U|^p)$ is the $p$th absolute moment of a
standard normal
variable $U$ (other versions are possible; one may, e.g., take
any product of $k$ increments, with powers adding up to $2$).
\end{longlist}

The first method has been introduced by Mancini in \cite{Manc1}, the second
one by Barndorff-Nielsen and Shephard in \cite{BS}. Both provide estimators
which converge in probability to $C_1$, upon
rather weak assumptions on the jumps.

The question of the rate of convergence, though, is still open, and we
quickly review the known results, in the case of truncated realized
volatility. One needs
the following assumption, where $r$ is a number in $[0,2]$:

{\renewcommand{\theass}{(L-$r$)}
\begin{ass}\label{assLr}
The variables $\sup_{t\leq1} |b_t|$, 
$\sup_{t\leq1} c_t$ and\break
$\sup_{t\leq1} \int(|x|^r\wedge1) F_t(dx)$ are almost surely
finite.
\end{ass}}%

The larger $r$ is, the weaker Assumption \ref{assLr} is.
(L-2) is a very weak assumption for an It\^o semimartingale, whereas
\ref{assLr} when $r<2$ puts restrictions on the jump activity,
and is slightly stronger than saying that the
Blumenthal--Getoor index of $X$ (or, jump activity index) is not bigger
than $r$. In particular, (L-$1$) is slightly stronger than the property
of the jumps to be summable on each finite interval, for example, the
jump part to
have trajectories of finite variation. Note that
a stable process of index $\beta\in(0,2)$ satisfies \ref{assLr} for all
$r>\beta$, but not for $r\leq\beta$.\vspace*{2pt}

When \ref{assLr} holds for some $r<1$, the estimators
$\wC(v_n)_1$ enjoy exactly the same CLT as in (\ref{3}) with
$\wC(v_n)$ in place of $[X,X]_t$, with the
same limit, provided we have
%
\begin{equation}
\label{6} v_n\asymp1/n^\varpi\qquad\mbox{with }
\frac{1}{4-2r}<\varpi<\frac{1}2.
\end{equation}
When \ref{assLr} holds for some $r\geq1$, the CLT with rate $\sqrt{n}$ no longer
holds for $\wC(v_n)$, but we have when $v_n\asymp1/n^\varpi$ with
$\varpi\in(0,1/r)$:
%
\begin{equation}
\label{7} 0<\varpi<\tfrac{1}2\qquad\Longrightarrow\quad n^{\varpi(2-r)} \bigl(
\wC(v_n)_1-C_1\bigr)\stackrel{\mathbb {P}} {
\longrightarrow}0
\end{equation}
(convergence in probability). These results are shown in~\cite{J08}, and
Mancini in~\cite{Manc2} has proved that when the
jumps of $X$ are those of a stable process with index $\beta$ [so
\mbox{\ref{assLr}} holds for all $r>\beta$, but not for $r=\beta$], and when
$\beta\geq1$,
the estimator converges exactly at rate $n^{\varpi(2-\beta)}$, in the sense
that the sequence $n^{\varpi(2-\beta)} (\wC(v_n)_1-C_1)$ converges
to a
nontrivial limit
(in probability, and not in law, in this case): this rate is less than
$\sqrt{n}$, as it is in (\ref{7}), and no proper CLT is available in
this case.

Turning now to multipowers, we have analogous results, except that
one needs stronger assumptions: basically, \ref{assLr} plus the fact that
the process $c_t$ is also an It\^o semimartingale, and never
vanishes: the CLT for $\wC(k,n)_1$ holds when $r<1$, with $\sqrt{2}$ replaced
by a suitable (bigger) constant depending on $k$; see \cite{BGJPS}.
When $r=1$, Vetter in \cite{Vetter-10} proves that there is a CLT at rate
$\sqrt{n}$ with a
nonvanishing bias term. When $r>1$ nothing is formally known, but the
presence of the bias term when $r=1$ suggests that for $r>1$ the rate
is less than $\sqrt{n}$.

\section{The results}\label{sec-R}

We are in a nonparametric setting, in which the process $X$ is not
specified [apart from the fact that it satisfies \ref{assLr} for some $r$],
and even the space $(\Omega,\mathcal{F},(\mathcal{F}_t)_{t\geq
0},\mathbb{P})$ is not specified. The meaning
of ``optimality'' or ``rate-optimality'' is not {a priori} clear;
and, to begin with, even the quantity to estimate, namely $C_1$, depends
of course on the space $(\Omega,\mathcal{F},(\mathcal{F}_t)_{t\geq
0},\mathbb{P})$ and on $X$.

A possible setting is as follows. We consider a family $\mathcal{S}$
of It\^o
semimartingales satisfying \ref{assLr}, each one being defined on its own
filtered space $(\Omega,\mathcal{F},(\mathcal{F}_t)_{t\geq
0},\mathbb{P})$, and the quantity to estimate is the associated
integrated volatility $C(X)_1$. Each $X$ in $\mathcal{S}$ takes its
values, as a process,
in the Skorokhod space~$\mathbb{D}^1$ of all c\`adl\`ag functions on
$\mathbb{R}_+$, and
the image by $X$ of the observed $\sigma$-field $\sigma
(X_{i/n}\dvtx i=0,\ldots,n)$
is the $\sigma$-field $\mathcal{D}_n=\sigma(x(i/n)\dvtx i=0,1,\ldots,n)$
of $\mathbb{D}^1$. For any
$X\in\mathcal{S}$ we denote by $\mathbb{P}^n_X$ the restriction to
$\mathcal{D}_n$ of the
law of $X$.

An\vspace*{2pt} estimator at stage $n$ is a $\mathcal{D}_n$-measurable function
$X\mapsto\wC(X)^n_i$ on $\mathbb{D}^1$.
We say that a sequence $\wC^n_1$ of such estimators achieves the
uniform rate
$w_n$ (with $w_n\to\infty$) on $\mathcal{S}$, for estimating $C_1$,
if the family
$w_n(\wC(X)^n_1-C(X)_1)$ is tight, uniformly in $n$ and in $X\in
\mathcal{S}$, that is,
$|\wC(X)^n_1-C(X)_1|=O_P(w_n^{-1})$ uniformly in $X\in\mathcal{S}$.

Of course, if $\mathcal{S}^r$ denotes the set of \textit{all} It\^o
semimartingales
satisfying \ref{assLr}, there cannot be any uniform rate because, to begin with,
the variables $C(X)_1$ are not uniformly tight when $X$ runs through
$\mathcal{S}^r$:
we need to restrict our attention to subfamilies of $\mathcal{S}^r$
which are
``bounded'' in some sense. In view of the formulation of~\ref{assLr}, it is
natural to consider, for any $A>0$, the class
%
\begin{eqnarray}\label{R-1} %
&& \mathcal{S}^r_A = \mbox{the
set of all It\^o semimartingales with}
\nonumber\\[-10pt]\\[-10pt]
&& \mbox{$|b_t|+c_t+\int(|x|^r\wedge1)
F_t(dx) \leq  A
\mbox{ for all $t$}$}.\nonumber
\end{eqnarray}
We also denote by $\mathcal{S}^{r,L}_A$ the subclass of all L\'evy processes
belonging to $\mathcal{S}^r_A$.

The main result of this paper is the following theorem.

%
\begin{theo}\label{TR-1} Let $r\in[0,2)$ and $A>0$. Any uniform rate $w_n$
for estimating $C(X)_1$, within the class $\mathcal{S}^{r,L}_A$, hence
also within
the bigger class $\mathcal{S}^r_A$, satisfies (up to a
multiplicative constant, of course)
%
\begin{equation}
\label{R-2} w_n\leq\rho_n:=\cases{ \sqrt{n}, &\quad
if $r\leq1$,
\vspace*{3pt}\cr
(n\log n)^{(2-r)/2},&\quad if $r>1$.}
\end{equation}
\end{theo}

The results recalled in the previous section show that the truncated
estimators $\wC(v_n)_1$ (which are estimators in the sense specified
above) achieve the rate $\rho_n$ when $r<1$, and at least
$n^{\varpi(2-r)}$ when $r\geq1$, for any $X$ satisfying \ref{assLr}.
We indeed have (slightly) more:

%
\begin{theo}\label{TR-2} Let $r\in[0,2)$ and $A>0$, and take
$v_n\asymp1/n^\varpi$.
The truncated estimators $\wC(v_n)_1$ have the uniform rate
$w_n$ below, within $\mathcal{S}^r_A$, for estimating~$C(X)_1$,
%
\begin{equation}
\label{R-3} w_n=\cases{ \sqrt{n}, &\quad if $r<1$ and $
\displaystyle\frac{1}{4-2r}\leq\varpi <\frac{1}2$,
\cr
n^{\varpi(2-r)},& \quad if $r\geq1$ and $0<\varpi<\dfrac{1}2$.}
\end{equation}
\end{theo}

When $r<1$, the truncated estimators $\wC(v_n)_1$ achieve the uniform rate
$\sqrt{n}$, and as seen in the previous section they even enjoy a CLT.
When $r\geq1$ we have the uniform rate $n^{\varpi(2-r)}$, although
for any
given $X$ we indeed have a ``faster'' rate, as seen in (\ref{7}); however,
this faster rate
is not uniform in $X\in\mathcal{S}^r_A$, as could be seen by taking a sequence
of L\'evy processes with characteristics $(0,1,G_n)$, with
$\int(|x|^r\wedge1) G_n(dx)\leq1$ (so $X^n\in\mathcal{S}^r_1$ for
all $n$),
but such that $\sup_n\int_{\{|x|\leq\varepsilon\}}|x|^r G_n(dx)$
does \textit{not} tend to $0$ as $\varepsilon\to0$.

We then conclude that the truncated estimators
are uniformly rate optimal when $r<1$, and otherwise they approach the bound
$\rho_n$, up to $n^{-\varepsilon}$ with $\varepsilon>0$ arbitrarily
small, upon
choosing $\varpi$ close enough to $\frac{1}2$.

Let us finally show that on the restricted class $\mathcal{S}^{r,L}_A$
of L\'evy
processes the rate $\rho_n$ of (\ref{R-2}) can be achieved exactly
and thus
constitutes the exact minimax optimal rate: this means that for any
$r\in[0,2)$
and any $A>0$ one can find estimators for $C(X)_1$ enjoying the uniform rate
$\rho_n$. When $r<1$, we already know this (even for the much
larger class $\mathcal{S}^r_A$) by the previous theorem, but for all
$r\in[0,2)$
we can
construct estimators with the uniform rate $\rho_n$ on $\mathcal
{S}^{r,L}_A$ as
follows. For any process~$X$, we
consider the empirical characteristic function of the
increments, at each stage $n$ (below $u\in\mathbb{R}$):
%
\begin{equation}
\label{R-10} \widehat{\phi}_n(u)=\frac{1}n\sum
_{j=1}^ne^{iu\Delta_j^nX}.
\end{equation}
Then we set
%
\begin{equation}
\label{R-11} \wC'(u)_1=-\frac{2n}{u^2} \bigl(\log\bigl|
\widehat{\phi}_n(u)\bigr| \bigr) 1_{\{\widehat{\phi}_n(u)\neq0\}}.
\end{equation}

%
\begin{theo}\label{TR-3} For all $A>0$ and $r\in[0,2)$, the estimators
$\wC'(u_n)_1$ with
%
\begin{equation}
\label{R-30} u_n=\cases{ \sqrt{n}, &\quad if $r\leq1$,
\vspace*{3pt}\cr
\sqrt{(r-1)n\log n} /\sqrt{A}, &\quad if $r>1$}
\end{equation}
attain the uniform rate $\rho_n$ for estimating $C(X)_1$, within
the class $\mathcal{S}^{r,L}_{A}$ of L\'evy processes.
\end{theo}

%
\begin{rem}\label{R1} When $r\leq1$ the estimators $\wC'(u_n)_1$ are likely
to enjoy a Central Limit theorem with rate $\rho_n$, and with a bias
when $r=1$.

When $r>1$, and upon examining the proof [see (\ref{C-13}) and (\ref{C-15}),
e.g.], the estimation error $\wC'(u_n)_1-C(X)_1$ is the sum\vspace*{1pt}
of a random part, which is easily seen to enjoy a CLT with rate
$n^{(2-r)/2}\log n$,
and\vspace*{1pt} a nonrandom part equal to $\Gamma_n=\frac{2\rho_n}{u_n^2}
\int(1-\cos(u_nx)) F(dx)$, where $F$ is the L\'evy measure of the L\'evy
process $X$ under consideration. It\vspace*{1pt} turns out that $|\rho_n\Gamma
_n|\leq
\int(u_n^{-r}\wedge|x|^r) F(dx)$ tends\vspace*{1pt} to $0$ by Lebesgue's theorem,
so, for any given $X$ we indeed
have \mbox{$\rho_n(\wC'(u_n)_1-C(X)_1)\to0$} in probability: this convergence
is of course not uniform in $X\in S^{r,L}_A$, otherwise the conclusion
of Theorem \ref{TR-1} would be violated. Now, depending on whether
$\rho_n\Gamma_n(\log n)^{r/2}$ converges or diverges---and\vspace*{1pt} both
occurrences are
possible---we have a CLT with rate $\rho_n(\log n)^{r/2}$, or we have
a slower effective rate (still at least~$\rho_n$, of course) with the
normalized error converging in probability to a nontrivial limit.

Note that the argumentation is in line with the standard nonparametric error
decomposition in a bias and variance part. Our estimator uses that the
characteristic exponent for high frequencies $u_n$ separates the Brownian
from the jump part according to the ratio $u_n^2/u_n^r$. We have reliable
empirical access to this exponent only up to frequency $u_n$ (otherwise the
stochastic error explodes due to a Gaussian deconvolution setting). So far,
we do not know whether this spectral approach yields the same optimal
rate on the larger class $\mathcal{S}^r_A$.
\end{rem}

\section{Proofs}\label{sec-P}

\subsection{Proof of Theorem \texorpdfstring{\protect\ref{TR-1}}{3.1}}

\textit{The bound} $w_n\leq\sqrt{n}$. For proving this bound, it is enough to
show that it already holds on the subclass $\mathcal{S}^{\mathrm{BM}}_A$ of all Brownian
motions with unit variance $c\leq A$ (so $\mathcal{S}^{\mathrm{BM}}_A\subset
\mathcal{S}^{r,L}_A$
for all $r\in[0,2]$).

In this case, and as already mentioned in the \hyperref[sec-Intro]{Introduction}, the increments
follow the parametric model $N(0,c/n)^{\otimes n}$ with parameter $c$ running
through $[0,A]$, for which
the LAN property holds with rate $\sqrt{n}$, and the result follows.

\textit{The bound} $w_n\leq(n\log n)^{(2-r)/2}$ \textit{when} $r\in(0,2)$.
By scaling, if the result holds for one $A>0$, it holds for all $A>0$.
Hence, in order to find a bound on the uniform rate $w_n$ on $\mathcal
{S}^{r,L}_A$,
hence {a fortiori} on $\mathcal{S}^r_A$, it
is enough to construct two sequences $X^n$ and $Y^n$ of L\'evy processes
belonging to $\mathcal{S}^{r,L}_K$ for $n\geq2$ and some constant
$K$, with the
following two properties, where $a_n=(n\log n)^{-(2-r)/2}$:
%
\begin{eqnarray}
&\bullet& \mbox{we have $C\bigl(X^n\bigr)_1=1+a_n$
and $C\bigl(Y^n\bigr)_1=1$ identically}, \label{P-1}
\\
&\bullet&\mbox{the total variation distance between $\mathbb
{P}^n_{X^n}$ and $\mathbb{P}^n_{Y^n}$
tends to $0$}. \label{P-2}
\end{eqnarray}

Indeed, letting $\wC(X)_1$ be a sequence of estimators with uniform rate
\mbox{$w_n\to\infty$} on\vspace*{2pt} $\mathcal{S}^r_A$ (or, even, on $\mathcal
{S}_A^{r,L}$), the two sequences
$w_n(\wC(X^n)^n_1-(1+a_n))$ and $w_n(\wC(Y^n)^n_t-1)$ are tight under
$\mathbb{P}^n_{X^n}$ and $\mathbb{P}^n_{Y^n}$, respectively, by (\ref
{P-1}). Then
(\ref{P-2})
implies that the sequence $w_n(\wC(Y^n)^n_1-(1+a_n))$ is also
tight under $\mathbb{P}^n_{Y^n}$. This is possible only if the sequence
$w_na_n$ is
bounded. So $1/a_n$ is an upper bound for any uniform rate on $\mathcal
{S}^{r,L}_K$
(up to a multiplicative constant, of course).

The proof of (\ref{P-1}) and (\ref{P-2}) is divided into several steps:
\begin{longlist}[(3)]
\item[(1)] We take L\'evy processes $X^n$ and $Y^n$
with respective characteristics $(0,1+a_n,F_n)$ and $(0,1,G_n)$, with
L\'evy measures $F_n,G_n$ satisfying
%
\begin{equation}
\label{C-1} \int\bigl(|x|^r\wedge1\bigr) F_n(dx)\leq
K,\qquad\int\bigl(|x|^r\wedge1\bigr) G_n(dx)\leq K
\end{equation}
for some constant $K$ (below constants change from line to line, and may
depend on $r$, and are all denoted as $K$).

By construction, we have (\ref{P-1}) and $X^n,Y^n\in\mathcal{S}^{r,L}_K$
for a constant $K$ [which may differ from the one in (\ref{C-1})], and we
need to choose the above
measures $F_n$ and $G_n$ in such a way that (\ref{P-2}) is satisfied.

\item[(2)] We take $u_n=2/a_n^{1/(2-r)}=2\sqrt{n\log n}$ and the even
functions $h_n\in C^2(\mathbb{R})$ defined for $u\geq0$ by
\[
h_n(u)=a_n \bigl(1_{\{u\leq u_n\}}+e^{-(u-u_n)^3}
1_{\{u>u_n\}} \bigr).
\]

We use the following convention for the Fourier transform, namely
$\mathcal{F}g(u)=\int e^{iux} g(x) \,dx$, so the inverse is
$\mathcal{F}^{-1}h(x)=\frac{1}{2\pi} \int e^{-iux} h(u) \,du$. We
also denote
as $f^{(q)}$ the $q$th derivative of any $q$-differentiable function $f$.

Since $h_n^{(q)}\in\mathbb{L}^p$ for all $p\geq1$ and $q=0,1,2$, we can
define $H_n=\mathcal{F}^{-1}h_n$, and we have $h_n^{(q)}=i^q\mathcal
{F}^{-1}H_{n,q}$, where
$H_{n,q}(x)=x^qH_n(x)$. By the Plancherel identity we deduce
%
\begin{eqnarray}\label{C-3}
&& \|H_{n}\|_{\mathbb{L}^2}\leq Ka_nu_n^{1/2}
\leq Ka_n^{(3-2r)/(4-2r)}, \qquad q=1,2
\nonumber\\[-8pt]\\[-8pt]
&&\qquad \Rightarrow\quad \|H_{n,q} \|_{\mathbb{L}^2}\leq\bigl\|h_{n}^{(q)}\bigr\|_{\mathbb{L}^2}\leq
Ka_n.\nonumber
\end{eqnarray}
Then\vspace*{-2pt} the Cauchy--Schwarz inequality applied to the functions
$\frac{1}{\sqrt{1+x^2}}$ and $H_n(x)\sqrt{1+x^2}$ yields
%
\begin{equation}
\label{C-4} \int\bigl|H_n(x)\bigr| \,dx\leq K\bigl(1+a_n^{(3-2r)/(4-2r)}
\bigr)<\infty
\end{equation}
[note that $\|H_n\|_{\mathbb{L}^1}$ is bounded in $n$ when $r\leq
3/2$, but not
otherwise; we also have $H_n(0)>a_nu_n\to\infty$]. Therefore, the two measures
\[
F_n(dx)=\frac{|H_n(x)|}{x^2} \,dx,\qquad G_n(dx)=F_n(dx)+
\frac{H_n(x)}{x^2} \,dx
\]
are nonnegative and integrate $x^2$, hence are L\'evy measures.\vadjust{\goodbreak}

This construction will satisfy (\ref{P-2}) mainly because the
definition of
the two L\'evy measures and the constant value of $h_n$ for $|u|\leq u_n$
imply that the difference between the two characteristic exponents vanishes
for $|u|\leq u_n$, as we shall prove next.

\item[(3)] Splitting the integration domain
into the sets $\{|u|\leq u_n\}$ and $\{|u|>u_n\}$ in the
integral $\int e^{-iux} h_n(u) \,du$, we get
\begin{eqnarray*}
\bigl|H_n(x)\bigr|&\leq& Ka_n \biggl(\frac{|\sin(u_nx)|}{|x|} +1 \biggr)
\\
&\leq& Ka_n \biggl(u_n 1_{\{|x|\leq1/u_n\}} +\frac{1}{|x|}
1_{\{1/u_n<|x|\leq1\}}+1_{\{|x|>1\}} \biggr).
\end{eqnarray*}
In\vspace*{2pt} turn, the integral $\int\frac{|x|^r\wedge1}{x^2} |H_n(x)| \,dx$
can be split into integrals on the sets $\{|x|\leq1/u_n\}$,
$\{1/u_n<|x|\leq1\}$ and $\{|x|>1\}$, and recalling $1<r<2$ we deduce from
the above that
\[
\int\frac{|x|^r\wedge1}{x^2} \bigl|H_n(x)\bigr| \,dx\leq Ka_n
\bigl(u_n^{2-r}+1\bigr)\leq K.
\]
It follows that the measures $F_n$ and $G_n$ satisfy (\ref{C-1}), and it
remains to prove~(\ref{P-2}).

\item[(4)] We\vspace*{2pt} denote by $\phi_n$ and $\psi_n$ the characteristic functions of
$X^n_{1/n}$ and $Y^n_{1/n}$, and $\eta_n=\phi_n-\psi_n$. These
functions are
real (because $H_n$ is an even function) and are given by
\begin{eqnarray*}
\phi_n(u)&=&\exp \biggl(-\frac{1}{2n} \bigl(u^2+a_nu^2+2
\Wphi_n(u) \bigr) \biggr),
\\
\psi_n(u)&=&\exp \biggl(-\frac{1}{2n} \bigl(u^2+2
\Wphi_n(u)+2\Weta _n(u) \bigr) \biggr),
\end{eqnarray*}
where
\begin{eqnarray*}
\Wphi_n(u)&=&\int\bigl(1-\cos(ux)\bigr)
\frac{|H_n(x)|}{x^2} \,dx,
\\
\Weta_n(u)&=&\int\bigl(1-\cos(ux)\bigr)
\frac{H_n(x)}{x^2} \,dx.
\end{eqnarray*}

We\vspace*{1pt} proceed to studying $\Wphi_n$ and $\Weta_n$. Equation~(\ref{C-3}) applied with
$q=1,2$ implies that $\Wphi_n$ and $\Weta_n$ are twice differentiable.
First, we have $\Wphi'_n(u)=\int\sin(ux) \frac{|H_n(x)|}{x} \,dx$, hence~(\ref{C-4}) yields
%
\begin{eqnarray}\label{C-5}
0&\leq&\Wphi_n(u)\leq K\bigl(1+a_n^{(3-2r)/(4-2r)}
\bigr) u^2,
\nonumber\\[-8pt]\\[-8pt]
\bigl|\Wphi'_n(u)\bigr|&\leq& K \bigl(1+a_n^{(3-2r)/(4-2r)}\bigr) |u|.\nonumber
\end{eqnarray}
Second, $\Weta''_n(u)=\int\cos(ux) H_n(x) \,dx=h_n(u)$,
whereas $\Weta(0)=\Weta'_n(0)=0$, and this yields
%
\begin{eqnarray}\label{C-6}
|u|\leq u_n &\quad\Rightarrow\quad &\Weta_n(u)= \frac{a_nu^2}2,\qquad \Weta'_n(u)=a_nu, \nonumber
\\[-8pt]
\\[-8pt]
|u|\geq u_n &\quad\Rightarrow\quad &\bigl|\Weta_n(u)\bigr|\leq \frac{a_nu^2}2,\qquad \bigl|\Weta'_n(u)\bigr|\leq
a_n|u|.
\nonumber
\end{eqnarray}

\item[(5)] Since $X^n$ and $Y^n$ have a nonvanishing Gaussian part,
the variables $X^n_{1/n}$ and $Y^n_{1/n}$ have\vspace*{1pt} densities, denoted by $f_n$
and $g_n$, and we set $k_n=f_n-g_n$. Since $X^n$ and $Y^n$
are L\'evy processes,
the variation distance between $\mathbb{P}^n_{X^n}$ and $\mathbb
{P}^n_{Y^n}$ is not more
than $n$ times $\int|k_n(x)| \,dx$, and we are thus left to show that
$n\int|k_n(x)| \,dx\to0$.

To check this, we use the same argument as for (\ref{C-4}): if $k_{n,1}(x)=
xk_n(x)$, by the Cauchy--Schwarz inequality we have $\int|k_n(x)| \,dx
\leq K(\|k_n\|_{\mathbb{L}^2}+\break \|k_{n,1}\|_{\mathbb{L}^2})$, whereas
$\eta_n=\mathcal{F}k_n$
and also, since $\eta_n$ is twice differentiable, $\eta'_n=i\mathcal
{F}k_{n,1}$.
By Plancherel identity, it is thus enough to prove that
%
\begin{equation}
\label{C-7} n^2\int\bigl|\eta_n(u)\bigr|^2 \,du\to 0,
\qquad n^2\int\bigl|\eta'_n(u)\bigr|^2 \,du
\to 0.
\end{equation}

We have $\Wphi_n\geq0$ and $\Wphi_n+\Weta_n\geq0$,
which implies $\phi_n(u)\leq e^{-u^2/2n}$ and $\psi_n(u)\leq e^{-u^2/2n}$,
whereas $2\Beta_n(u)=a_nu^2$ if $|u|\leq u_n$ and
$|2\Beta_n(u)|\leq a_nu^2$ if $|u|>u_n$ by~(\ref{C-6}). Therefore,
\begin{eqnarray*}
\bigl|\eta_n(u)\bigr| &=& \psi_n(u) \biggl|1-\frac{\phi_n(u)}{\psi_n(u)} \biggr|
\\
&=& \psi_n(u) \bigl|1-e^{-(a_nu^2-2\Weta_n(u))/(2n)} \bigr| \leq \frac{a_nu^2}{2n}
e^{-u^2/2n} 1_{\{|u|>u_n\}},
\end{eqnarray*}
and also, upon using (\ref{C-5}),
\begin{eqnarray*}
\bigl|\eta'_n(u)\bigr|&=&\frac{1}n \bigl|
\bigl(u+ua_n+\Wphi'_n(u)\bigr)
\phi_n(u) -\bigl(u+\Wphi'_n(u)+
\Beta'_n(u)\bigr)\psi_n(u) \bigr|
\\
&\leq&\frac{1}n \bigl(a_n|u|e^{-u^2/2n}+\bigl|
\Weta'_n(u)\bigr|e^{-u^2/2n} +\bigl|u+\Wphi'_n(u)\bigr|
\bigl|\eta_n(u)\bigr| \bigr) 1_{\{|u|>u_n\}}
\\
&\leq&Ka_n \frac{|u|}n e^{-u^2/2n} \biggl(1+
\bigl(1+a_n^{({3-2r})/({4-2r})}\bigr)\frac{u^2}n \biggr)
1_{\{|u|>u_n\}}.
\end{eqnarray*}
Now, since $u_n=2\sqrt{n\log n}$, we have
$\int_{\{|u|>u_n\}} (\frac{u^2}{n} )^q e^{-u^2/n} \,du\leq
K\frac{(\log n)^{q-1}}{n^{7/2}}$ for $q=1,2,3$. Since
further $a_n^{(3-2r)/(4-2r)}/\sqrt{n}\to0$, we deduce
\[
\int\bigl|\eta_n(u)\bigr|^2 \,du\leq K\frac{\log n}{n^{7/2}}, \qquad
\int\bigl|\eta'_n(u)\bigr|^2 \,du\leq K
\frac{(\log n)^2}{n^{7-1/2}}.
\]
Then (\ref{C-7}) follows, and the proof is complete.
\end{longlist}

\subsection{Proof of Theorem \texorpdfstring{\protect\ref{TR-2}}{3.2}}

The proof requires several steps:
\begin{longlist}[(3)]
\item[(1)] Any $X\in\mathcal{S}^r_A$ can be written as follows, on some space
$(\Omega,\mathcal{F},(\mathcal{F}_t)_{t\geq0},\mathbb{P})$:
%
\begin{eqnarray}
\label{C-10} X_t&=&X_0+\int_0^tb_s
\,ds+\int_0^t\sqrt{c_s}
\,dW_{s}\nonumber
\\
&&{} +\int_0^t\!\int_E \delta(s,z) 1_{\{\|\delta(s,z)\|\leq1\}} (\mu-\nu) (ds,dz)
\\
&&{} +\int_0^t\!\int_E
\delta(s,z) 1_{\{\|\delta(s,z)\|>1\}} \mu(ds,dz).
\nonumber
\end{eqnarray}
Here, $b$ and $c$ are as in \ref{assLr}, and $W$ is a standard Brownian motion,
and $\mu$ is a Poisson random measure on $\mathbb{R}_+\times\mathbb
{R}$ with intensity
measure $\nu(dt,dz)=dt\otimes dz$, and $\delta=\delta(\omega,t,z)$ is
a predictable function on $\Omega\times\mathbb{R}_+\times\mathbb
{R}$. The connection between
$\delta$ and $F_t$ is that $F_{\omega,t}$ is the image of Lebesgue
measure by the
map $z\mapsto\delta(\omega,t,z)$, restricted to $\mathbb
{R}\setminus\{0\}$.

We use the decomposition $X=X'+Y+Z$, where
\[
X'_t=X_0+\int_0^tb_s
\,ds+\int_0^t\sqrt{c_s}
\,dW_{s}
\]
and $Y$ and $Z$ are, respectively, the last two terms in (\ref{C-10}).

With $w_n$ given by (\ref{R-3}), it is clearly enough to prove that,
for some constant $K$ only depending on $A,r,\varpi$ (as will be all
constants $K$ below, changing from line to line), we have
%
\begin{equation}
\label{C-11} \mathbb{E} \bigl(\bigl|\wC(v_n)_1-C_1\bigr|
\bigr)\leq K/w_n.
\end{equation}

\item[(2)] Here, we recall estimates on the increments of $X'$ and $Y$, the later
coming from Lemmas 2.1.5 and 2.1.6 of \cite{JP}, and where $p>0$ is arbitrary
(the constants $K_p$ below depend on $p$ in addition to $r,A$). Namely, since
$\int_{\{|x|\leq1\}}|x|^r F_t(dx)\leq A$, we have uniformly in
$s\in[(i-1)/n,i/n]$:
%
\begin{eqnarray}
\label{C-20}
\mathbb{E}\bigl(\bigl|X'_s-X'_{(i-1)/n}\bigr|^p
\bigr)&\leq& K_pn^{-p/2},
\nonumber\\[-8pt]\\[-8pt]
\mathbb{E}\bigl(\bigl|Y_s-Y_{(i-1)/n}\bigr|^p
\bigr)&\leq& Kn^{-(p/r)\wedge1}.\nonumber
\end{eqnarray}
We will also use the following estimates, which follow from the property
$F_t(\{x\dvtx |x|>1\})\leq A$ and from the fact that if $\Delta^n_iZ\neq0$ there
is at
least one jump of $Z$ within the interval $ (\frac{i-1}n,\frac{i}n ]$
(this estimate follows from Lemma 2.1.7 of \cite{JP} applied
to the counting process $\sum_{s\leq t}1_{\{\Delta Z_s\neq0\}}$):
%
\begin{equation}
\label{C-21} \mathbb{P}\bigl(\Delta^n_iZ\neq0\bigr)
\leq\frac{K}n.
\end{equation}

\item[(3)] With the notation (\ref{2}), It\^o's formula yields
$[X',X']^n_1-C_1=U_n+V_n$, where
\begin{eqnarray*}
U_n&=&\sum_{i=1}^n\mathbb{E}
\bigl(\zeta^n_i\mid\mathcal {F}_{(i-1)/n}\bigr),
\\
\zeta^n_i&=&2\int_{(i-1)/n}^{i/n}
\bigl(X'_s-X'_{(i-1)/n}
\bigr)b_s \,ds,
\\
V_n&=&\sum_{i=1}^n
\xi^n_i,
\\
\xi^n_i&=&2\int
_{(i-1)/n}^{i/n}\bigl(X'_s-X'_{(i-1)/n}
\bigr)\sqrt{c_s} \,dW_s +\zeta^n_i-
\mathbb{E}\bigl(\zeta^n_i\mid\mathcal{F}_{(i-1)/n}
\bigr).
\end{eqnarray*}
Equation~(\ref{C-20}) yields
\[
\bigl|\mathbb{E}\bigl(\zeta^n_i\mid\mathcal{F}_{(i-1)/n}
\bigr) \bigr|\leq K/n^{3/2},\qquad \mathbb{E}\bigl(\bigl(\xi^n_i
\bigr)^2\bigr)+\mathbb{E}\bigl(\bigl(\zeta^n_i
\bigr)^2\bigr)\leq K/n^2,
\]
whereas $\mathbb{E}(\xi^n_i\mid\mathcal{F}_{(i-1)/n})=0$. Thus we
have $\mathbb{E}(|U_n|)\leq
K/\sqrt{n}$
and $\mathbb{E}(V_n^2)\leq K/n$, implying
%
\begin{equation}
\label{C-12} \mathbb{E} \bigl(\bigl|\bigl[X',X'
\bigr]^n_1-C_1\bigr| \bigr)\leq K/\sqrt{n}.
\end{equation}
Therefore, it remains to prove that
%
\begin{equation}
\label{C-110} \mathbb{E} \bigl(\bigl|\wC(v_n)_1-
\bigl[X',X'\bigr]^n_1\bigr| \bigr)
\leq K/w_n.
\end{equation}

\item[(4)] Consider the case $r<1$ first. By Lemma 13.2.6 of \cite{JP}, applied
with $k=1$ and $F(x)=x^2$, hence $s'=2$ and $m=s=p'=1$ and $\theta=0$
(with the notation of this lemma), we have
\[
\mathbb{E} \Biggl( \Biggl|\wC(v_n)_1-\sum
_{i=1}^n\bigl(\Delta^n_iX'
\bigr)^2 1_{\{
|\Delta^n_iX'|\leq v_n\}} \Biggr| \Biggr)\leq\frac{K}{n^{(2-r)\varpi}}\leq
\frac{K}{\sqrt{n}},
\]
where the last inequality follows from $\varpi\geq\frac{1}{4-2r}$. On
the other
hand, (\ref{C-20}) and Markov inequality yield
$\mathbb{E}((\Delta^n_iX')^2 1_{\{|\Delta^n_iX'|>v_n\}})\leq
K_p/n^{1+p(1-2\varpi)/2}$ for
any \mbox{$p>0$}, and upon taking $p=\frac{1}{1-2\varpi}$ we obtain
\[
\mathbb{E} \Biggl( \Biggl|\bigl[X',X'\bigr]^n_1-
\sum_{i=1}^n\bigl(\Delta^n_iX'
\bigr)^2 1_{\{
|\Delta^n_iX'|\leq v_n\}} \Biggr| \Biggr)\leq\frac{K}{\sqrt{n}}.
\]
These two estimates readily give (\ref{C-110}).

\item[(5)]  Now we turn to the case $r\geq1$. One has
$\wC(v_n)_1-[X',X']^n_1
=\break\sum_{j=1}^3U(j)_n$, where $U(j)_n=\sum_{i=1}^n\eta(j)^n_i$ and
\begin{eqnarray*}
\eta(1)^n_i&=&\bigl(\Delta^n_iX
\bigr)^2 1_{\{|\Delta^n_iX|\leq v_n\}}-\bigl(\Delta ^n_iX'
\bigr)^2-2\Delta^n_iX'
\Delta^n_iY,
\\
\eta(2)^n_i&=&2\mathbb{E}\bigl(\Delta^n_iX'
\Delta^n_iY\mid\mathcal {F}_{(i-1)/n}\bigr),\qquad
\eta(3)^n_i=2\Delta^n_iX'
\Delta^n_iY-\eta(2)^n_i.
\end{eqnarray*}
It\^o's formula yields, with the notation $\gamma_s=\int_{\{|z|\leq
1\}}
z^2 F_s(dz)$, and taking advantage of the facts that $Y$ and
$\int_0^t\sqrt{c_s} \,dW_s$ are two orthogonal martingales and that
$Y^2_t-\int_0^t\gamma_s \,ds$ is a martingale:
\begin{eqnarray*}
&& \eta(2)^n_i =2\mathbb{E} \biggl(\int
_{(i-1)/n}^{i/n}\bigl(X'_s-X'_{(i-1)/n}
\bigr)b_s \,ds \Big|\mathcal{F}_{(i-1)/n} \biggr)
\\
&& \mathbb{E} \bigl(\bigl(\Delta^n_iX'
\Delta^n_iY\bigr)^2\mid\mathcal
{F}_{(i-1)/n} \bigr)
\\
&&\qquad = \mathbb{E} \biggl(\int_{(i-1)/n}^{i/n}(Y_s-Y_{(i-1)/n})^2
c_s \,ds \Big|\mathcal{F}_{(i-1)/n} \biggr)
\\
&&\quad\qquad{} +2\mathbb{E} \biggl(\int_{(i-1)/n}^{i/n}
\bigl(X'_s-X'_{(i-1)/n}\bigr)
(Y_s-Y_{(i-1)/n})^2 b_s \,ds\Big|
\mathcal{F}_{(i-1)/n} \biggr)
\\
&&\quad\qquad{} +\mathbb{E} \biggl(\int_{(i-1)/n}^{i/n}
\bigl(X'_s-X'_{(i-1)/n}
\bigr)^2 \gamma_s \,ds \Big|\mathcal{F}_{(i-1)/n}
\biggr).
\end{eqnarray*}
Then standard estimates and (\ref{C-20}), plus H\"older's inequality, yield
(the first bound is a.s.)
\[
\bigl|\eta(2)^n_i\bigr|\leq\frac{K}{n^{3/2}},\qquad \mathbb{E}
\bigl(\bigl(\eta(3)^n_i\bigr)^2 \bigr)\leq
\frac{K}{n^2}.
\]
Since $\mathbb{E}(\eta(3)^n_i\mid\mathcal{F}_{(i-1)/n})=0$, these
estimates yield
$|U(2)_n|\leq K/\sqrt{n}$ and\break  $\mathbb{E}(U(3)^2_n)\leq K/n$, hence
it is enough
to show that $\mathbb{E}(|U(1)_n|)\leq K/w_n$.

\item[(6)] Recalling $r\geq1$, the following inequality is easy to check, for
$x,y,z\in\mathbb{R}$ and $v\in(0,1/4]$:
\begin{eqnarray*}
&& \bigl|(x+y+z)^2 1_{\{\bigr|x+y+z|\leq v\}}-x^2-2xy |
\\
&&\qquad \leq
2v^2 1_{\{z\neq0\}}+6|xy| 1_{\{|x|>v/2\}} +6x^2
1_{\{|x|>v/2\}}+16 v^{2-r}|y|^r.
\end{eqnarray*}
It follows that $|\eta(1)^n_i|\leq K\sum_{j=1}^5\xi(j)^n_i$, where
\begin{eqnarray*}
\xi(1)^n_i&=&v_n^2
1_{\{\Delta^n_iZ\neq0\}},\qquad \xi(2)^n_i=\bigl|
\Delta^n_iX'\Delta^n_iY\bigr|
1_{\{|\Delta^n_iX'|>v_n/2\}},
\\
\xi(3)^n_i&=&\bigl(\Delta^n_iX'
\bigr)^2 1_{\{|\Delta^n_iX'|\geq v_n/2\}},\qquad \xi(4)^n_i=v_n^{2-r}
\bigl|\Delta^n_iY\bigr|^r.
\end{eqnarray*}
Equation~(\ref{C-21}) yields $\mathbb{E}(\xi(1)^n_i)\leq K/n^{1+2\varpi}$,
and (\ref{C-20})
yields $\mathbb{E}(\xi(4)^n_i)\leq K/n^{1+(2-r)\varpi}$. Another\vspace*{2pt}
application of
(\ref{C-20}), plus H\"older and Markov inequalities, give us
$\mathbb{E}(\xi(j)^n_i)\leq K_{p}/n^{1+p(1-2\varpi)/2}$ for $j=2,3$.
Upon taking
$p$ large enough, we obtain
\[
\mathbb{E}\bigl(\xi(j)^n_i\bigr)\leq K/nw_n
\]
for $j=1,2,3,4,5$. We deduce $\mathbb{E}(|U(1)_n|)\leq K/w_n$, and the proof
is complete.
\end{longlist}

\subsection{Proof of Theorem \texorpdfstring{\protect\ref{TR-3}}{3.3}}

We\vspace*{1pt} let $X\in\mathcal{S}^{r,L}_A$, where $r\in[0,2)$ and $A>0$ are
given. The
characteristic triple of $X$ is $(b,c,F)$ and the characteristic function
of $X_{1/n}$ is
\[
\phi_n(u)=\exp \biggl(\frac{1}n \biggl(iub-\frac{cu^2}2+\int
\bigl(e^{iux}-1- iux 1_{\{|x|\leq1\}}\bigr) F(dx) \biggr) \biggr).
\]
Then $|\phi_n(u_n)|=e^{(-{1}/({2n})) (cu_n^2+\gamma_n )}$, where
$\gamma_n=2\int(1-\cos(u_nx)) F(dx)$. As soon as $n$ is large
enough we have
$u_n\geq1$, hence, since $1-\cos y\leq1\wedge
y^2\leq|y|^r\wedge1$,
\begin{eqnarray*}
0&\leq&\gamma_n\leq2\int\bigl(|u_nx|^r\wedge1
\bigr) F(dx) \leq2u_n^r\int\bigl(|x|^r\wedge1
\bigr) F(dx)
\\
& \leq&2u_n^2\int\bigl(|x|^r\wedge1
\bigr) F(dx).
\end{eqnarray*}
Because $c+\int(|x|^r\wedge1) F(dx)\leq A$ by hypothesis, and in
view of
the form of $u_n$ in~(\ref{R-30}), by singling out the two cases
$r\leq1$
and $r>1$ this implies that, with $\Gamma=e^A$,
%
\begin{equation}
\label{C-13} \frac{1}{|\phi_n(u_n)|}=e^{u_n^2(c+\gamma_n)/2n} \leq\Gamma n^{(r-1)^+/2}.
\end{equation}

The estimation error $\wC'(u_n)_1-c$ is the sum $G_n+H_n$ of the deterministic
and stochastic errors:
\begin{eqnarray*}
G_n&=& -\frac{2n}{u_n^2} \log\bigl|\phi_n(u_n)\bigr|-c=
\frac{\gamma_n}{u_n^2},
\\
H_n&=&\frac{2n}{u_n^2} \bigl(\log\bigl|\phi_n(u_n)\bigr|-
\bigl(\log\bigl|\widehat {\phi}_n(u_n)\bigr| \bigr)
1_{\{\widehat{\phi}_n(u_n)\neq0\}} \bigr).
\end{eqnarray*}
The previous estimates on $\gamma_n$ readily yield
%
\begin{equation}
\label{C-14} |G_n|\leq\frac{2A}{u_n^{2-r}}.
\end{equation}

Second, we study $H_n$. The variables $\exp(iu_n\Delta^n_jX)$ are i.i.d.
as $j$ varies, with\vspace*{1pt} modulus $1$ and expectation $\phi_n(u_n)$, hence
$V_n=\widehat{\phi}_n(u_n)-\phi_n(u_n)$ satisfies $\mathbb
{E}(|V_n|^2)\leq1/n$. In view
of (\ref{C-13}), on the set $\{|V_n|\leq1/n^{r/4}\}$
we have $|V_n/\phi_n(u_n)|\leq1/2$ and $\widehat{\phi}_n(u_n)=
V_n+\phi_n(u_n)\neq0$ as soon\vspace*{1pt} as $n\geq n_0=\break  (2\Gamma
)^{4/((2-r)\wedge
r)}$, in
which case we deduce, for some universal constant $K$:
\[
|H_n|=\frac{2n}{u_n^2} \biggl|\log \biggr|1+\frac{V_n}{\phi
_n(u_n)} \biggl| \biggr| \leq K
\frac{n|V_n|}{u_n^2 |\phi_n(u_n)|}.
\]
Henceforth, if $n\geq n_0$,
%
\begin{equation}
\label{C-15} \mathbb{E} \bigl(|H_n| 1_{\{|V_n|\leq1/n^{r/4}\}} \bigr)\leq \cases{
\displaystyle\frac{K\Gamma}{\sqrt{n}}, &\quad if $r\leq1$,
\vspace*{3pt}\cr
\displaystyle \frac{KA\Gamma}{ (r-1) n^{(2-r)/2} \log n},
&\quad if $r>1$.}
\end{equation}

Putting\vspace*{1pt} together (\ref{C-14}) and (\ref{C-15}), plus the fact
that $\mathbb{P}(|V_n|>1/n^{r/4})\leq1/n^{(2-r)/2}$ (by Bienaym\'
e--Tchebycheff
inequality)\vspace*{1pt} tends to zero, and the equality
$\wC'(u_n)_1-c=G_n+H_n$, we deduce that $\rho_n(\wC'(u_n)_1-c)$
[with the
notation~(\ref{R-2})] is tight, uniformly in $X\in\mathcal{S}^{r,L}_A$.




\printaddresses
\end{document}